# An embedding for the Kesten–Spitzer random walk in random scenery

**Endre Csáki**[*]

*Mathematical Institute of the Hungarian Academy of Sciences, Budapest, P.O.B. 127, H-1364, Hungary.* E-mail: `csaki@math-inst.hu`

**Wolfgang König**

*Technische Universität Berlin, Fachbereich Mathematik, Str. des 17 Juni 135, D-10623 Berlin, Germany.* E-mail: `koenig@math.tu-berlin.de`

**Zhan Shi**

*Laboratoire de Probabilités, Université Paris VI, 4 Place Jussieu, F-75252 Paris Cedex 05, France.* E-mail: `shi@ccr.jussieu.fr`

**Summary.** For one-dimensional simple random walk in a general i.i.d. scenery and its limiting process, we construct a coupling with explicit rate of approximation, extending a recent result for Gaussian sceneries due to Khoshnevisan and Lewis. Furthermore, we explicitly identify the constant in the law of iterated logarithm.

**Keywords.** Local time, random walk in random scenery, Brownian motion in Brownian scenery, strong approximation.

**1991 Mathematics Subject Classification.** 60J55; 60J15; 60J65.

[*] Research supported by the Hungarian National Foundation for Scientific Research, Grant No. T 016384 and T 019346.

# 1. Introduction

Let $\sigma = \{\sigma_x\}_{x \in \mathbb{Z}}$ (sometimes also written $\{\sigma(x)\}_{x \in \mathbb{Z}}$) denote a sequence of independent and identically distributed real-valued random variables such that

$$(1.1) \qquad \mathbb{E}(\sigma_0) = 0, \quad \mathbb{E}(\sigma_0^2) = 1 \quad \text{and} \quad \mathbb{E}(|\sigma_0|^p) < \infty, \quad \text{for all } p > 0.$$

Any realization of the sequence $\{\sigma_x\}_{x \in \mathbb{Z}}$ is called a "scenery". Let $S = \{S_k\}_{k \in \mathbb{N}_0}$ be a simple symmetric random walk on $\mathbb{Z}$ starting at $S_0 = 0$, independent of $\sigma$. The process $K = \{K(n)\}_{n \in \mathbb{N}_0}$, defined by

$$(1.2) \qquad K(n) = \sum_{k=0}^{n} \sigma(S_k), \qquad n \in \mathbb{N}_0,$$

is usually referred to as the **Kesten–Spitzer random walk in random scenery**, see [8] for more details. For example, the model can be viewed as follows: if a random walker has to pay the amount of $\sigma_x$ dollars whenever he visits the site $x$, then $K(n)$ is the total amount he pays during the first $n$ steps.

There is a continuous analogue for $K$ introduced and analyzed by Kesten and Spitzer [8]. To describe this, let $B = \{B(t); t \geq 0\}$ and $W = \{W(x); x \in \mathbb{R}\}$ be independent real-valued standard Brownian motions with $B(0) = W(0) = 0$. Let $\{L(t,x); t \geq 0, x \in \mathbb{R}\}$ denote the jointly continuous version of the local time process of $B$, in the sense that for any non-negative Borel function $f$,

$$\int_0^t f(B(s)) \, \mathrm{d}s = \int_{\mathbb{R}} f(x) L(t,x) \, \mathrm{d}x, \qquad t > 0,$$

see Trotter [16]. Now, define the process $G$, which will be called **Brownian motion in Brownian scenery**, by:

$$(1.3) \qquad G(t) = \int_{\mathbb{R}} L(t,x) \, \mathrm{d}W(x), \qquad t > 0.$$

It is proved by Kesten and Spitzer [8] that

$$(1.4) \qquad \left\{ n^{-3/4} K(\lfloor nt \rfloor); 0 \leq t \leq 1 \right\} \xrightarrow{\text{law}} \left\{ G(t); 0 \leq t \leq 1 \right\},$$

where "$\xrightarrow{\text{law}}$" stands for weak convergence in law (in some functional space; for example in the space of bounded functions on $[0,1]$ endowed with the uniform topology).

Clearly, the process $G$ is self-similar in the sense that for any $a > 0$,

$$(1.5) \qquad G(a \cdot) \stackrel{\text{law}}{=} a^{3/4} G(\cdot),$$

with "$\stackrel{\text{law}}{=}$" denoting identity in distribution.

It is natural to ask whether (1.4) holds in a stronger sense. For example, is it possible to obtain a strong approximation for $K$ by $G$? This problem was studied by Khoshnevisan and Lewis [9] in the special case where the random scenery $\sigma$ is Gaussian. More precisely, they proved that, if $\sigma_0$ is a Gaussian $\mathcal{N}(0,1)$ variable, then (possibly in an enlarged probability space) one can construct a pair of $K$ and $G$ such that, with probability one, for any $\varepsilon > 0$,

$$(1.6) \qquad \max_{0 \leq m \leq n} |K(m) - G(m)| = o(n^{1/2+\varepsilon}), \qquad n \to \infty.$$

In view of the self-similarity in (1.5), we immediately recover (1.4) from (1.6), in the case of Gaussian scenery.

It is one main aim of this paper to extend (1.6) to any random scenery satisfying (1.1). The precise formulation is as follows.

**Theorem 1.1.** *Let $\{\sigma_x\}_{x \in \mathbb{Z}}$ be a random scenery satisfying (1.1). Possibly in an enlarged probability space, there exists a coupling for $K$ and $G$, such that, with probability one, for any $\varepsilon > 0$,*

$$\max_{0 \leq m \leq n} |K(m) - G(m)| = o(n^{5/8+\varepsilon}), \qquad n \to \infty.$$

**Remark.** The rate $o(n^{5/8+\varepsilon})$ in Theorem 1.1 for a general Gaussian scenery is not as good as the one in (1.6) for a Gaussian scenery. This originates from the fact that from Brownian motion, it is easier to construct an embedding for Gaussian variables than for arbitrary variables.

A straightforward consequence of Theorem 1.1 is a law of the iterated logarithm (LIL) for $K$. Indeed, it is proved by Khoshnevisan and Lewis [9] that for some absolute constant $c_0 \in (0, \infty)$,

$$(1.7) \qquad \limsup_{t \to \infty} \frac{G(t)}{(t \log \log t)^{3/4}} = c_0, \qquad \text{a.s.}$$

Therefore, an application of Theorem 1.1 yields

$$(1.8) \qquad \limsup_{n \to \infty} \frac{K(n)}{(n \log \log n)^{3/4}} = c_0, \qquad \text{a.s.}$$

Our second main result identifies the exact value of the constant $c_0$.

**Theorem 1.2.** *Both* (1.7) *and* (1.8) *hold with* $c_0 = 2^{5/4}/3$.

For other properties of $K$ and/or $G$, we refer to Khoshnevisan and Lewis [9] and the references therein. We also mention the recent work of Khoshnevisan and Lewis [10] where some interesting open problems are raised, with partial answers in Xiao [17].

The rest of the paper is organized as follows. The proof of Theorem 1.1 is carried out in Section 2 in two steps. In fact, we need two embeddings of different natures. The first embedding consists of constructing a random walk $S$ from the Brownian motion $B$, whereas the second is a construction of a random scenery $\sigma$ from the Brownian scenery $W$. In order to clarify the embeddings, we first outline the method in the beginning of Section 2 by formulating the two main steps. There we also explain the relation of our work to Khoshnevisan and Lewis [9].

Theorem 1.2 is proved in Section 3.

# 2. Proof of Theorem 1.1

## 2.1. OUTLINE

Recall that $\sigma$ denotes the random scenery satisfying (1.1), $S$ the random walk and $B$ the Brownian motion (having local times $L$) and $W$ the Brownian scenery.

Introduce the number of the walker's visits to $x$ until time $n$,

$$(2.1) \qquad \xi(n,x) = \sum_{k=0}^{n} \mathbb{1}_{\{S_k = x\}}, \qquad n \in \mathbb{N}_0,\, x \in \mathbb{Z},$$

which is often referred to as the local time of the random walk $S$. Then (1.2) can be rewritten as

$$(2.2) \qquad K(n) = \sum_{x \in \mathbb{Z}} \sigma_x\, \xi(n,x), \qquad n \in \mathbb{N}_0.$$

We shall prove Theorem 1.1 in two steps: First, we approximate $\xi(n,x)$ by $L(n,x)$ (the walker's imbedding), and secondly, we approximate $\sigma_x$ by $dW(x)$ (the scenery's imbedding). The precise formulation is as follows.

**Proposition 2.1.** *There is a coupling of $\sigma$, $S$ and $B$ such that $\sigma$ is independent of $(S, B)$ and such that, with probability one, for any $\varepsilon > 0$,*

$$(2.3) \qquad \sum_{x \in \mathbb{Z}} \sigma_x \left(\xi(n, x) - L(n, x)\right) = o(n^{1/2+\varepsilon}), \qquad n \to \infty.$$

**Proposition 2.2.** *There is a coupling of $\sigma$, $S$, $B$ and $W$ such that $(\sigma, W)$ is independent of $(S, B)$ and such that, with probability one, for any $\varepsilon > 0$,*

$$(2.4) \qquad \sum_{x \in \mathbb{Z}} \sigma_x L(n, x) - \int_{\mathbb{R}} L(n, x) \, \mathrm{d}W(x) = o(n^{5/8+\varepsilon}), \qquad n \to \infty.$$

It is straightforward to see that Theorem 1.1 follows from these two propositions and their proofs. In the next two subsections we shall prove the propositions.

We say a few words about the difference between our approach and the one adopted by Khoshnevisan and Lewis [9].

In the walker's embedding (proof of Proposition 2.1), the main difficulty consists of estimating the moments of a certain random variable in order to apply the Borel-Cantelli lemma. If the random scenery $\{\sigma_x\}_{x \in \mathbb{Z}}$ is Gaussian, then the random variable in question is conditionally Gaussian (given the random walk), which allowed Khoshnevisan and Lewis [9] to obtain accurate estimates for the moments. In the general case, the Gaussian techniques break down, but we succeed in the moments estimate by means of a general inequality for moments, taken from Petrov [12].

If the scenery is Gaussian, then it is itself equal (in distribution) to the increments of the Brownian scenery $W$ at integer times, so there was no need for a scenery's embedding in [9]. In the case of a general scenery, we shall use the Skorokhod embedding to approximate the scenery by the increments of $W$ at certain random times. This will leave us with a remainder which finally leads to the term $n^{5/8+\varepsilon}$ in Theorem 1.1 (instead of $n^{1/2+\varepsilon}$ in the case of a Gaussian scenery).

## 2.2. Proof of Proposition 2.1: The walker's embedding

According to a theorem by Révész [13] (see also Chapter 10 of Révész [14]), one can construct a random walk $S$ from the Brownian motion $B$ such that, with probability one, for every $\varepsilon > 0$,

$$(2.5) \qquad \sup_{x \in \mathbb{Z}} |\xi(n, x) - L(n, x)| = o(n^{1/4+\varepsilon}), \qquad n \to \infty,$$

where $\xi$ and $L$ denote the local times of $S$ resp. $B$. In particular, we may assume that $\sigma$ and $(S, B)$ are independent.

Define

$$(2.6) \qquad I(N, n) = \sum_{x=-N}^{N} \sigma_x(\xi(n, x) - L(n, x)), \qquad N, n \in \mathbb{N}.$$

To finish the proof via the first Borel-Cantelli lemma, it is sufficient to show that, for every $\varepsilon > 0$, the probabilities $\mathbb{P}(|I(\infty, n)| > n^{1/2+\varepsilon})$ are summable over $n \in \mathbb{N}$ (with an obvious definition of $I(\infty, n)$).

Pick $\varepsilon > 0$. First note that, according to the classical LIL's for Brownian motion and simple random walk,

$$(2.7) \qquad \limsup_{t \to \infty} \frac{\sup_{0 \le s \le t} |B(s)|}{(2t \log \log t)^{1/2}} = 1 = \limsup_{n \to \infty} \frac{\max_{0 \le k \le n} |S_k|}{(2n \log \log n)^{1/2}}, \qquad \text{a.s.}$$

Hence, we have, with probability one, $L(n, x) = 0 = \xi(n, x)$ for all $|x| \ge n^{1/2+\varepsilon}$ and sufficiently large $n$, and therefore $I(\infty, n) = I(\lfloor n^{1/2+\varepsilon} \rfloor, n)$. Hence, it suffices to show the summability of $\mathbb{P}(|I(\lfloor n^{1/2+\varepsilon} \rfloor, n)| > n^{1/2+\varepsilon})$.

Pick some $p \ge 2$. We are going to use (see Petrov [12], p. 62) that, with some constant $c_1 = c_1(p) > 0$, for any sequence $\{X_i\}_{i \in \mathbb{N}}$ of independent (but not necessarily identically distributed) mean-zero random variables, we have the estimate

$$(2.8) \qquad \mathbb{E}\left( \left| \sum_{i=1}^{N} X_i \right|^p \right) \le c_1 N^{p/2-1} \sum_{i=1}^{N} \mathbb{E}(|X_i|^p), \qquad N \in \mathbb{N}.$$

We apply this fact to the variables $\sigma_x(\xi(n, x) - L(n, x))$, $x = -N, \ldots, N$, conditioned on $(S, B)$, and obtain, using the independence of $\sigma$ and $(S, B)$, that

$$(2.9) \qquad \mathbb{E}(|I(N, n)|^p) \le c_1 (2N + 1)^{p/2-1} \sum_{x=-N}^{N} \mathbb{E}(|\sigma_0|^p) \mathbb{E}(|\xi(n, x) - L(n, x)|^p).$$

It is also proved by Khoshnevisan and Lewis [9] that, for the construction of the walk $S$ from the motion $B$ we are using, the distance of their local times is also small in $L^p$ sense, more precisely, there is a constant $c_2 = c_2(p) > 0$ such that

$$(2.10) \qquad \sup_{x \in \mathbb{Z}} \mathbb{E}(|\xi(n, x) - L(n, x)|^p) \le c_2 \, n^{p/4}, \qquad n \in \mathbb{N}.$$

Using this on the r.h.s of (2.9), we obtain that

$$(2.11) \qquad \mathbb{E}(|I(N, n)|^p) \le \mathcal{O}(N^{p/2}) \mathcal{O}(n^{p/4}), \qquad n, N \to \infty.$$

Now we use the Chebyshev inequality and apply (2.11) to $N = \lfloor n^{1/2+\varepsilon} \rfloor$ to get

$$\mathbb{P}\left(\left|I(\lfloor n^{1/2+\varepsilon}\rfloor, n)\right| > n^{1/2+\varepsilon}\right) \leq n^{-p/2-p\varepsilon} \mathbb{E}\left(\left|I(\lfloor n^{1/2+\varepsilon}\rfloor, n)\right|^p\right)$$
(2.12)
$$= \mathcal{O}(n^{-p\varepsilon/2}), \qquad n \to \infty,$$

which is summable for $p > 2/\varepsilon$. This ends the proof of Proposition 2.1.

## 2.3. Proof of Proposition 2.2: The scenery's embedding

Let $B = \{B(t); t \geq 0\}$ and $W = \{W(x); x \in \mathbb{R}\}$ be a Brownian motion resp. scenery satisfying $B(0) = 0 = W(0)$. Let $\sigma = \{\sigma_x\}_{x \in \mathbb{Z}}$ be a random scenery satisfying (1.1). In the proof of Proposition 2.1 we constructed a simple random walk $S$ from $B$ (thus, independent of $W$) whose local times satisfy (2.3).

We are using now the classical Skorokhod embedding (see Breiman [1, Theorem 13.8]) which ensures the existence of i.i.d. non-negative random variables $T_i$, $i \in \mathbb{Z}$, with $\mathbb{E}(T_i) = \mathbb{E}(\sigma_i^2) = 1$ such that

$$\left\{W\left(\sum_{i=1}^n T_i\right)\right\}_{n \in \mathbb{Z}} \stackrel{\text{law}}{=} \left\{\sum_{x=1}^n \sigma_x\right\}_{n \in \mathbb{Z}},$$
(2.13)

with the notation $\sum_{i=1}^0 a_i \stackrel{\text{def}}{=} 0$ and $\sum_{i=1}^n a_i \stackrel{\text{def}}{=} a_{-1} + \cdots + a_n$ for negative $n$. Since $\sigma_0$ possesses all moments, also the variables $T_i$ do. For brevity, we write

$$\varrho(n) = \varrho_n \stackrel{\text{def}}{=} \sum_{i=1}^n T_i, \qquad n \in \mathbb{Z},$$
(2.14)

so that we have

$$\{\widetilde{\sigma}_n\}_{n \in \mathbb{Z}} \stackrel{\text{def}}{=} \{W(\varrho_n) - W(\varrho_{n-1})\}_{n \in \mathbb{Z}} \stackrel{\text{law}}{=} \{\sigma_n\}_{n \in \mathbb{Z}}.$$
(2.15)

Note that we have constructed $\widetilde{\sigma} = (\widetilde{\sigma}_x)_{x \in \mathbb{Z}}$ from $W$ and may therefore assume that $(\widetilde{\sigma}, W)$ and $(S, B)$ are independent.

For $N \in \mathbb{N}$, abbreviate

$$J(N, n) \stackrel{\text{def}}{=} \int_0^{\varrho(N)} L(n, x) \, dW(x) - \sum_{j=1}^N \widetilde{\sigma}_j L(n, j), \qquad n \in \mathbb{N}.$$
(2.16)

To finish the proof it is sufficient to show that, with probability one, for any $\varepsilon > 0$ (using an obvious notation),

$$J(\infty, n) = o(n^{5/8+\varepsilon}), \qquad n \to \infty.$$
(2.17)

since the integral over negative $x$ and the sum over negative $j$ are handled in the same way. Note that

(2.18) $$J(N,n) = \int_0^{\varrho(N)} A_n(x)\,dW(x), \qquad N, n \in \mathbb{N},$$

where

$$A_n(x) \stackrel{\text{def}}{=} L(n,x) - L(n,j), \qquad \text{if } x \in (\varrho_{j-1}, \varrho_j],$$

for $j \in \mathbb{N}$. By the Dambis–Dubins–Schwarz representation theorem for continuous local martingales (see for example Theorem V.1.6 of Revuz and Yor [15]), there exists, for every $n \in \mathbb{N}$, a Brownian motion $\{\mathbb{W}_n(t);\, t \geq 0\}$ such that

(2.19) $$\int_0^t A_n(x)\,dW(x) = \mathbb{W}_n\left(\int_0^t A_n^2(x)\,dx\right), \qquad t \geq 0.$$

In particular, we have, with probability one,

(2.20) $$J(N,n) = \mathbb{W}_n\left(\int_0^{\varrho(N)} A_n^2(x)\,dx\right), \qquad N, n \in \mathbb{N}.$$

On the other hand, using a well-known estimate for the Gaussian tail gives that

$$\mathbb{P}\left(\max_{1 \leq i \leq n} \sup_{0 \leq s \leq t} |\mathbb{W}_i(s)| > \alpha\right) \leq n\mathbb{P}\left(\sup_{0 \leq s \leq t} |W(s)| > \alpha\right) \leq 4n\exp\left(-\frac{\alpha^2}{2t}\right),$$

which, combined with an application of the first Borel–Cantelli lemma and monotonicity, yields that for any $a > 0$, with probability one,

(2.21) $$\max_{1 \leq i \leq n} \sup_{0 \leq s \leq n^a} |\mathbb{W}_i(s)| = \mathcal{O}\left(n^{a/2}(\log n)^{1/2}\right), \qquad n \to \infty.$$

(We mention that it is possible to obtain an estimate more accurate than (2.21), by means of Theorem 1.2 of Deheuvels and Révész [4]).

We are going to apply a result from Csáki et al. [2] which says that for any $a \geq 0$ and $\varepsilon > 0$, with probability one,

(2.22) $$\sup_{|x-y| \leq t^a} |L(t,x) - L(t,y)| = o(t^{1/4+a/2+\varepsilon}), \qquad t \to \infty.$$

Furthermore, since the sequence $(T_i)_{i \in \mathbb{N}}$ is i.i.d. with all moments finite and $\mathbb{E}(T_1) = 1$, the classical Hartman–Wintner LIL implies that, with probability one,

(2.23) $$\varrho(n) = n + \mathcal{O}\left((n \log \log n)^{1/2}\right), \qquad n \to \infty.$$

Now fix $b > 1/2$ and $\varepsilon > 0$. For any $0 < x \le \varrho(\lfloor n^b \rfloor)$, there exists $1 \le j \le \lfloor n^b \rfloor$ such that $x \in (\varrho_{j-1}, \varrho_j]$. Since by (2.23), $|x - j| \le n^{b/2+\varepsilon}$ and $x \le 2n^b$, we can apply (2.22) to conclude that, with probability one,

$$(2.24) \qquad \sup_{0 \le x \le \varrho(\lfloor n^b \rfloor)} |A_n(x)| = o(n^{(1+b)/4+\varepsilon}), \qquad n \to \infty.$$

Therefore, by (2.23),

$$(2.25) \qquad \int_0^{\varrho(\lfloor n^b \rfloor)} A_n^2(x)\,dx \le \varrho(\lfloor n^b \rfloor) \sup_{0 \le x \le \varrho(\lfloor n^b \rfloor)} A_n^2(x) = o(n^{(1+3b)/2+2\varepsilon}), \qquad n \to \infty.$$

Going back to (2.20), and by means of (2.21),

$$(2.26) \qquad J(\lfloor n^b \rfloor, n) = \mathcal{O}\big(n^{(1+3b)/4+\varepsilon}(\log n)^{1/2}\big) = o\big(n^{(1+3b)/4+2\varepsilon}\big).$$

Since $b > 1/2$, the usual LIL (recalled in (2.7)) implies that, for all large $n$, $L(n, x) = L(n, j) = 0$ for $x \ge \varrho(\lfloor n^b \rfloor)$ and $j \ge \lfloor n^b \rfloor$ and hence $J(\infty, n) = J(\lfloor n^b \rfloor, n)$. We have thus proved that

$$(2.27) \qquad J(\infty, n) = o\big(n^{(1+3b)/4+2\varepsilon}\big), \qquad n \to \infty,$$

for any $b > 1/2$. Since $(1 + 3b)/4$ can be made as close to $5/8$ as possible, we have proved that (2.17) holds for every $\varepsilon > 0$. Therefore, the proof of Proposition 2.2 is complete.

## 3. Proof of Theorem 1.2

### 3.1. OUTLINE

In this section, we describe how to get the exact value of the constant $c_0$ in (1.7) and (1.8), which was left open by Khoshnevisan and Lewis [9]. Let $L(t, x)$ as before denote Brownian local time. Khoshnevisan and Lewis [9] have proved that the value of $c_0$ is determined by

$$c_0 = \frac{2}{(27\,\zeta)^{1/4}},$$

where the constant $\zeta \in (0, \infty)$ is defined in terms of the Brownian self-intersection local time,

$$X_t \stackrel{\text{def}}{=} \int_{\mathbb{R}} L(t, x)^2\,dx, \qquad t > 0,$$

as follows:

$$(3.1) \qquad \zeta \stackrel{\text{def}}{=} -\lim_{\lambda \to \infty} \frac{1}{\lambda^2} \log \mathbb{P}(X_1 > \lambda);$$

see (6.1), (5.14) and Corollary 5.6 in [9].

Thus, Theorem 1.2 is proved as soon we have proved that $\zeta = 3/2$. In Subsection 3.2 we give an analytic proof for this fact, using Mansmann's [11] large deviation result for the exponential moments of $X_t$.

In Subsection 3.3 we show how some tools from stochastic analysis can be used to prove at least the inequality $\zeta \leq 3/2$. Unfortunately, we have not been able also to derive the opposite inequality $\zeta \geq 3/2$ (which is the much harder one) by similar means.

**Remark.** We note that our identification of $\zeta$ also can be used to identify the constant in the LIL for the process $(X_t)_{t>0}$ in Remark 1.2.1. in Csörgő et al [3], more precisely, we have there

$$\limsup_{t \to \infty} \frac{X_t}{(t^3 \log \log t)^{1/2}} = \left(\frac{2}{3}\right)^{1/2}, \qquad \text{a.s.}$$

## 3.2. ANALYTICAL PROOF FOR $\zeta = 3/2$

A particular case of Kasahara's [7] Tauberian theorem says that (3.1) is equivalent to

$$(3.2) \qquad \lim_{a \to \infty} \frac{1}{a^2} \log \mathbb{E}(e^{aX_1}) = \frac{1}{4\zeta}.$$

(This can also be easily checked by adapting the proof of Cramér's theorem). By the Brownian scaling property, we have $X_1 \stackrel{\text{law}}{=} t^{-3/2} X_t$, and therefore (3.2) is equivalent to

$$(3.3) \qquad \lim_{t \to \infty} \frac{1}{t} \log \mathbb{E}\left(\exp\left\{\frac{2}{t} X_t\right\}\right) = \frac{1}{\zeta}.$$

But the l.h.s. of (3.3) has been investigated by Mansmann [11] in his study of the polaron problem. Based on a general large deviation principle due to Donsker and Varadhan [5], Mansmann proved that

$$(3.4) \qquad \text{l.h.s. of (3.3)} = \sup_{\varphi}\left(2 \int_{\mathbb{R}} \varphi^4(x) \, dx - \frac{1}{2} \int_{\mathbb{R}} \varphi'(x)^2 \, dx\right),$$

where the supremum is taken over all absolutely continuous functions $\varphi \colon \mathbb{R} \to \mathbb{R}$ such that $\int_{\mathbb{R}} \varphi^2(x) \, dx = 1$. It is also proved in [11] that the maximizer on the right hand side of (3.4) is given by

$$\varphi_*(x) = \frac{1}{\cosh(2x)}, \qquad x \in \mathbb{R}.$$

(In pp. 94, 119 and 124 of [11], it was claimed that $\varphi_*(x) = 2^{-1/2}/\cosh(2x)$, but an inspection on Lemma 6.9 of the same reference reveals that the correct choice for $\varphi_*$ is $1/\cosh(2x)$). Since $4\varphi_*^3 + \frac{1}{2}\varphi_*'' = 2\varphi_*$ and since $\int_{\mathbb{R}} \varphi_*^4(x)\, dx = 2/3$, the right hand side of (3.4) equals 2/3. Combining (3.3) and (3.4), we arrive at our assertion $\zeta = 3/2$.

### 3.3. PROBABILISTIC PROOF FOR $\zeta \leq 3/2$

We are going to prove the inequality $\zeta \leq 3/2$, using identifications of the laws of some stochastic processes constructed from Brownian motion.

Let $\tau \stackrel{\text{def}}{=} \sup\{t < 1 : B(t) = 0\}$, the last passage time at 0 before 1 of the Brownian motion $B$. It is well-known that the process

$$\Lambda(t) \stackrel{\text{def}}{=} \frac{B(t\tau)}{\sqrt{\tau}}, \qquad t \in [0,1],$$

is a standard Brownian bridge, independent of $\tau$. Let $L_\Lambda(t,x)$ denote the local time process of $\Lambda$. A straightforward application of the occupation times formula yields that

$$L_\Lambda(t,x) = \frac{1}{\sqrt{\tau}} L(t\tau, x\sqrt{\tau}), \qquad t \in [0,1],\, x \in \mathbb{R}.$$

Therefore,

$$X_\tau = \int_{\mathbb{R}} L(\tau,x)^2\, dx = \tau^{3/2} \int_{\mathbb{R}} L_\Lambda(1,x)^2\, dx.$$

Now recall that

$$\lim_{\lambda \to \infty} \frac{1}{\lambda^2} \log \mathbb{P}\left( \int_{\mathbb{R}} L_\Lambda(1,x)^2\, dx > \lambda \right) = -\frac{3}{2}.$$

This was proved in Csörgő et al. [3], by means of Jeulin's [6] characterization of the local time of the normalized Brownian excursion. Since $\mathbb{P}(\tau > 1 - \varepsilon) > 0$ for any $\varepsilon > 0$, and since $\tau$ is independent of $L_\Lambda$, we conclude that

(3.5) $$\lim_{\lambda \to \infty} \frac{1}{\lambda^2} \log \mathbb{P}(X_\tau > \lambda) = -\frac{3}{2}.$$

Note that $X_\tau \leq X_1$ because $\tau \leq 1$, hence (3.5) yields the assertion $\zeta \leq 3/2$.

Unfortunately, we have not been able to use this approach for deriving the opposite inequality, $\zeta \geq 3/2$. It is intuitively clear that, on the event that $X_1$ is very large, $B(1)$ should be very close to zero and $\tau$ should be very close to 1. We have not been able to turn this idea into an honest proof. But this part of the proof is anyway the much harder one, as is seen from an inspection of Section 5 in [11].

# Acknowledgements

W. K. gratefully acknowledges the hospitality of The Fields Institute Toronto where part of the work has been carried out. Cooperation between E. Cs. and Z. S. was supported by the joint French-Hungarian Intergovernmental Grant 'Balaton' (grant no. F25/97). The authors also wish to acknowledge the support of the Paul Erdős Summer Research Center of Mathematics, Budapest.

# References

[1] Breiman, L.: *Probability.* Addison-Wesley, Reading, 1968.

[2] Csáki, E., Csörgő, M., Földes, A. and Révész, P.: The local time of iterated Brownian motion. *J. Theoretical Probab.* 9 (1996) 717–743.

[3] Csörgő, M., Shi, Z. and Yor, M.: Some asymptotic properties of the local time of the uniform empirical process. *Bernoulli* (to appear)

[4] Deheuvels, P. and Révész, P.: On the coverage of Strassen-type sets by sequences of Wiener processes. *J. Theoretical Probab.* 6 (1993) 427–449.

[5] Donsker, M.D. and Varadhan, S.R.S.: Asymptotics for the polaron. *Comm. Pure Appl. Math.* 36 (1983) 505–528.

[6] Jeulin, T: Applications du grossissement de filtrations à l'étude des temps locaux du mouvement brownien. In: *Grossissements de Filtrations: Exemples et Applications* (T. Jeulin and M. Yor, eds.). Lecture Notes in Mathematics 1118 pp. 197–304. Springer, Berlin, 1985.

[7] Kasahara, Y.: Tauberian theorems of exponential type. *J. Math. Kyoto Univ.* 18 (1978) 209–219.

[8] Kesten, H. and Spitzer, F.: A limit theorem related to a new class of self similar processes. *Z. Wahrsch. Verw. Gebiete* 50 (1979) 5–25.

[9] Khoshnevisan, D. and Lewis, T.M.: A law of the iterated logarithm for stable processes in random scenery. *Stoch. Proc. Appl.* 74 (1998) 89–121.

[10] Khoshnevisan, D. and Lewis, T.M.: Iterated Brownian motion and its intrinsic skeletal structure. (preprint)

[11] Mansmann, U.: The free energy of the Dirac polaron, an explicit solution. *Stochastics Stochastics Rep.* 34 (1991) 93–125.

[12] Petrov, V.V.: *Limit Theorems of Probability Theory.* Oxford Science Publications, Oxford, 1995.

[13] Révész, P.: Local time and invariance. In: *Analytical Methods in Probability Theory* (D. Dugué et al., eds.). Lecture Notes in Mathematics 861 pp. 128–145. Springer, Berlin, 1981.

[14] Révész, P.: *Random Walk in Random and Non–Random Environments.* World Scientific, Singapore, 1990.

[15] Revuz, D. and Yor, M.: *Continuous Martingales and Brownian Motion.* Springer, Berlin, 2nd ed., 1994.

[16] Trotter, H.F.: A property of Brownian motion paths. *Illinois J. Math.* 2 (1958) 425–433.

[17] Xiao, Y.: The Hausdorff dimension of the level sets of stable processes in random scenery. (preprint)